\documentclass[a4paper,final,12pt]{article}

\usepackage{amsthm}
\usepackage{graphicx,amsmath,amssymb,url}
\usepackage[usenames]{color}
\usepackage[citecolor=blue]{hyperref}
\hypersetup{
    colorlinks=true,
    linkcolor=blue,
    urlcolor=cyan,
}
\usepackage[width=15cm]{geometry}

\newtheorem{proposition}{Proposition}[section]
\newtheorem{theorem}[proposition]{Theorem}
\newtheorem{corollary}[proposition]{Corollary}
\newtheorem{lemma}[proposition]{Lemma}
\newtheorem{remark}[proposition]{Remark}

\newcommand{\wt}[1]{\widetilde{#1}}
\newcommand{\dt}[1]{{\bf #1}}
\newcommand{\ignore}[1]{}	% argument is ignored

\newcommand{\mt}{\rightarrow}

\newcommand{\beq}{\begin{equation}}
\newcommand{\eeq}{\end{equation}}
\newcommand{\beql}[1]{\begin{equation}\label{eq:#1}}
\newcommand{\eeql}{\end{equation}}
\newcommand{\beqarrays}{\begin{eqnarray*}}
\newcommand{\eeqarrays}{\end{eqnarray*}}

\newcommand{\blem}{\begin{lemma}}
\newcommand{\elem}{\end{lemma}}
\newcommand{\bleml}[1]{\begin{lemma} \label{lem:#1}}
\newcommand{\eleml}{\end{lemma}}
\newcommand{\blemT}[2]{\begin{lemma}[#1] \label{lem:#2}}
\newcommand{\elemT}{\end{lemma}}
\newcommand{\bthm}{\begin{theorem}}
\newcommand{\ethm}{\end{theorem}}
\newcommand{\bthml}[1]{\begin{theorem} \label{thm:#1}}
\newcommand{\ethml}{\end{theorem}}
\newcommand{\bthmT}[2]{\begin{theorem}[#1] \label{thm:#2}}
\newcommand{\ethmT}{\end{theorem}}
\newcommand{\brem}{\begin{remark}}
\newcommand{\erem}{\end{remark}}
\newcommand{\breml}[1]{\begin{remark} \label{rem:#1}}
\newcommand{\ereml}{\end{remark}}
\newcommand{\bcor}{\begin{corollary}}
\newcommand{\ecor}{\end{corollary}}
\newcommand{\bcorl}[1]{\begin{corollary} \label{cor:#1}}
\newcommand{\ecorl}{\end{corollary}}
\newcommand{\bpropl}[1]{\begin{proposition} \label{pro:#1}}
\newcommand{\epropl}{\end{proposition}}

\newcommand{\bpf}{\begin{proof}}
\newcommand{\epf}{\end{proof}\hfill \qed}

\newcommand{\refeq}[1]{(\protect\ref{eq:#1})}

\newcommand{\refLem}[1]{Lemma~\protect\ref{lem:#1}}
\newcommand{\refThm}[1]{Theorem~\protect\ref{thm:#1}}
\newcommand{\refFig}[1]{Figure~\protect\ref{fig:#1}}

\newcommand{\refSec}[1]{Section~\protect\ref{sec:#1}}

\newcommand{\refPro}[1]{Proposition~\protect\ref{pro:#1}}
\newcommand{\CC}{{\mathbb C}}

\newcommand{\ZZ}{{\mathbb Z}}

\newcommand{\as}{\textcolor{red}{\mathop{\mbox{\rm :=}}}}
\newcommand{\dd}{ ,\ldots , }

\newcommand{\ibp}{\subset }		% "proper \ib"
\newcommand{\abs}[1]{\left\lvert#1\right\rvert}       % magnitude
\newcommand{\set}[1]{\left\{ #1 \right\}}

\newcommand{\paren}[1]{\left( { #1 }\right)}	% parenthesis
\newcommand{\su}{\cup}

\newcommand{\sm}{\setminus}

	\newenvironment{prog}{\begin{tabbing}
	 xxxx\=xxxx\=xxxx\=xxxx\=xxxx\=xxxx\=xxxx\=xxxx\=xxxx\=xxxx\=xxxx\=xxxx\=xxxx\=
	\kill\\}{
		\end{tabbing}}

	\newenvironment{progb}[2][2]{ % Sep'2014: optional arg for headspace!
		\begin{center}
		  \fbox{\begin{minipage}{0.75\textwidth}
			\begin{prog}#2\end{prog}
	          \end{minipage}}
	        \end{center}
		}{}
	\newcommand{\vfigpdf}[3]{
            \begin{figure}[htb]
              \centering
              \includegraphics[scale=#3]{#2}
              \caption{#1}
              \label{fig:#2}
            \end{figure}
          }

\newcommand{\dia}{\text{dia}}
\newcommand{\diaG}{\text{dia}(G)}

\newcommand{\eitheta}{e^{i \theta}}

\newcommand{\fur}{F_{u,r}}

\newcommand{\stirtwo}[2]{\genfrac{\{}{\}}{0pt}{}{#1}{#2}}   %Stirling numbers of second kind
  \title{On The Roots of Independence Polynomial: Quantifying The Gap}
  \author{Om Prakash\\The Institute of Mathematical Sciences, HBNI, Chennai, India.\\
  \url{omprakashphd@outlook.com}\\
  Vikram Sharma\\The Institute of Mathematical Sciences, HBNI, Chennai, India.\\
  \url{vikram@imsc.res.in}}
\begin{document}
\maketitle

\begin{abstract}
The independence polynomial of a graph $G$ is the generating polynomial
corresponding to its independent sets of different sizes. More formally,
if $a_k(G)$ denotes the number of independent sets of $G$ of size $k$ then
     \[I(G,z) \as \sum_{k}^{} (-1)^k a_k(G) z^k.\] 
The study of evaluating $I(G,z)$ has several
deep connections to problems in combinatorics, complexity theory and statistical physics.
Consequently, the roots of the independence polynomial have been studied in detail.
In particular, many works have provided regions in the complex plane that are devoid
of any roots of the polynomial. One of the first such results showed a lower
bound on the absolute value of the smallest root $\beta(G)$ of the polynomial.
Furthermore, when $G$ is connected, Goldwurm and Santini 
established that $\beta(G)$ is a simple real root of $I(G,z)$ smaller than one.
An alternative proof was given by Csikv\'ari. Both proofs
do not provide a gap from $\beta(G)$ to the  smallest absolute value
amongst all the other roots of $I(G,z)$. In this paper, we quantify this gap.
\end{abstract}
  % \begin{keyword}
  %   Independence Polynomial, Root separation, Zero-free regions.
  % \end{keyword}

\section{Introduction}
\label{sec:intro}
Let $G=(V,E)$ be a simple undirected graph, that is, without loops and multiple edges, with
$V$ representing its set of $n$ vertices and $E$ its set of edges. An \dt{independent set}
of $G$ is a subset of vertices of $V$ such that there is no edge between any 
pair of vertices in the subset. Let $a_k(G)$ denote the number of independent sets of size
$k$ in $G$, where $a_0(G) \as 1$. 
Following the convention in \cite{csikvari:note:12}, we define the \dt{independence polynomial}
    \beql{indpoly}
    I(G,z) \as \sum_{k=0}^n (-1)^k a_k(G) z^k.
    \eeql
Let $\beta(G)$ denote the smallest real root of $I(G,z)$. 
It is well known that
such a root exists and is indeed in the interval $(0,1]$. Moreover, it is also known
that any other root $\rho$ of $I(G,z)$ is strictly greater than $\beta(G)$ in absolute value
\cite{goldwurm,csikvari:note:12}. We next describe their proofs in brief, but for this we need 
to introduce some notation. 

For a vertex $v \in V$, let $N_G(v)$ denote the set of neighbors of $v$ and $d(v)$
its degree in $G$; we will also use $N(v)$ if $G$ is clear from the context. 
The set of \dt{closed neighbors} of $v$ is $N[v] \as N(v) \su \set{v}$. Given a subset $S \ibp V$,
the graph $G\sm S$ denotes the subgraph of $G$ induced by the vertices $V\sm S$.
A key recursive property of the independence polynomial is the following:
For all vertices $u\in V$, we have
     \beql{recur}
     I(G, z) = I(G\sm u,z) - zI(G\sm N[u],z).
     \eeql

One of the key tools used in both \cite{goldwurm} and \cite{csikvari:note:12} is the Taylor series expansion
of $1/I(G,z)$ around the origin. In \cite{goldwurm}, the vertices of the graph are treated
as symbols of an alphabet and the edge relations as congruence relations on the alphabet, that is,
if two vertices have an edge then the corresponding alphabets can be swapped in any string.
The congruence relations impose an equivalence relation, called the trace monoid, 
on the set of all finite strings over the alphabet. It is well known that the $n$th 
coefficient of the power series $1/I(G,z)$ is the number of traces of length $n$ in the monoid.
Using the properties of the trace monoid, it is shown that 
the power series can be
expressed as a rational function where both the numerator and denominator correspond to
the characteristic polynomial of two positive matrices. Furthermore, the positive matrix 
that appears in the denominator dominates the one in the numerator entry wise. Therefore, its
largest eigenvalue is unique (due to Perron-Frobenius) and it does not appear as an eigenvalue of the numerator.
This establishes the uniqueness of the pole, which is also $\beta(G)$. 
We do not see an immediate way to quantify the proof.

In \cite{csikvari:note:12} it is shown,
using \refeq{indpoly}, that the coefficients
of the power series $1/I(G,z)$ are all positive. Moreover, from Pringsheim's theorem
we know that the radius of convergence of this power series is $\beta(G)$.
Now consider the power series $1/I(H,z)$ for any {\em proper subgraph} $H$ of $G$.
We can express
     \[\frac{1}{I(G,z)} = \frac{I(H,z)}{I(G,z)}\frac{1}{I(H,z)}.\]
By repeated applications of \refeq{indpoly}, along with induction, one can argue that
the coefficients of both $I(H,z)/I(G,z)$ and $1/I(H,z)$ are positive. 
Therefore, the coefficients of the series on the left-hand side above are greater 
than the coefficients of $1/I(H,z)$. Hence, $\beta(G) < \beta(H)$. 
By an inductive argument, it is then shown that $\beta(G)$ is in fact
a simple root of $I(G,z)$, when $G$ is connected. 
To show that any other root $\rho$ of $I(G,z)$ is strictly
greater than $\beta(G)$ in absolute value, they consider the function
     \[
     f_u(z) \as \frac{z I(G\sm N[u],z)}{I(G\sm u ,z)},
     \]
where $u \in V$. An inductive argument, similar to the one used for $1/I(G,z)$, shows
that the coefficients (except the constant coefficient)
of this power series are also positive. Now if $\rho$ is any other root of $I(G,z)$
with absolute value $\beta(G)$ %, which is smaller than the radius of convergence $\beta(G\sm u)$ of $f_u(z)$, 
then $f_u(\rho) = f_u(|\rho|) = 1$, that is, $f_u$ is periodic on the circle
$\beta(G) \eitheta$, where $\theta \in [0, 2\pi)$ and $i \as \sqrt{-1}$. 
Using the ``Daffodil Lemma''
from complex analysis \cite[p.~266]{fs}, this implies that the coefficients of $f_u(z)$ are a subset
of an arithmetic progression. However, since the $k$th coefficient of $f_u(z)$ asymptotically
is of the form $(\beta(G\sm u))^{-k}$, it does not satisfy an arithmetic progression 
and this gives us a contradiction that $\rho$
can have the same absolute value as $\beta(G)$. As the proof is by contradiction, it also
fails to quantify the gap between $\beta(G)$ and $\rho$.

The insight in this paper builds on the observation that since $f_u(z)$ is holomorphic 
for all $z$ such that $|z| < \beta(G\sm u)$, from the Maximum Modulus Principle \cite[p.~545]{fs}
we know that the largest absolute value of $f_u(z)$ on the disc $D(0, \beta(G)\eitheta)$ 
is attained on its boundary. Moreover, as the coefficients of $f_u(z)$ are positive,
a simple calculation shows that the maximum absolute value is attained on the
positive real axis, namely at $\beta(G)$. This gives an alternate proof to the ones given
above. To quantify this argument, we proceed in two steps: first, 
we show that the function is univalent in a neighborhood of $\beta(G)$, and second,
by constructing a disc around all points with absolute value $\beta(G)$, except in a certain
neighborhood of $\beta(G)$, where $f_u(z)$ does not take the value one. 
The two steps should intuitively hold, since in the
case of the former, as the derivative $f_u'(z)$ does not vanish at $\beta(G)$ there must
be a neighborhood of $\beta(G)$ where the function is injective; the latter case 
follows from the continuity of $f_u(z)$. The main challenge is to quantify these two intuitive ideas.
For this purpose, we need tools, such as, Smale's
$\gamma$-function  to study the function locally, and some simple results from 
complex analysis on radius of univalence of a function, such as $f_u(z)$. The main result
of the paper is the following:
\bthml{main}
Let $G$ be a connected graph on $n$ vertices. Then  the disc centered at the origin
\[D\paren{0, \beta(G) + \paren{ \frac{\beta(G)}{n}}^{O(n)}}\] 
contains only the smallest root $\beta(G)$ of $I(G,z)$.
\ethml

% Besides quantifying the gap, another motivation to study this ratio is whether it could
% help in our understanding of evaluating $I(G,z)$ in the complex plane. In the general case,
% there are results that tightly relate the absence of zeros of the independence polynomial
% to its polynomial time evaluation algorithms in certain regions of the complex plane \cite{patel-regts,srivastava}.
% These algorithms simultaneously study all graphs of fixed degree. But this does not
% explain that evaluating certain independence polynomials, such as those associated with 
% a cycle, a path or a complete bipartite graph can be efficiently evaluated in almost
% the entire complex plane. 
% In some of these cases, the gap between $\beta(G)$ and the absolute
% value of the next root is substantially larger than one, as we demonstrate with some examples. 
% Our larger aim is to understand the
% algorithmic implication of this gap in evaluating the independence polynomial.

In the next section, we describe the main results of the paper with intuitive details.
The proofs of these results are developed in the subsequent sections. 
The necessary preliminary results and definitions
from graph theory and complex analysis that are needed are provided in the appendix \refSec{prelim}.

\section{Main Results}
\label{sec:main-results}

Our focus will be to understand the properties of the following function: For any $u \in V$, let
     \beql{fuz}
     f_u(z) \as \frac{z I(G\sm N[u],z)}{I(G\sm u ,z)},     
     \eeql
For example, consider the star graph $S_n$ with one central vertex of degree $n$ connected to $n$ leaves.
It is not hard to verify that its independence polynomial is $(1-z)^n-z$. Now, if $u$ is the
``center'' vertex in $S_n$ then $f_u(z) = z/(1-z)^n$, but if $u$ is one of the leaves in $S_n$ then
     \beql{fuex}
     f_u(z) = z(1-z)^{n-1}/((1-z)^{n-1}-z).
     \eeql
If $N[u] = \set{u, u_1 \dd u_k}$, where $k=d(u)$, then 
    \[f_u(z) = z \frac{I(G\sm \set{u, u_1},z)}{I(G\sm u ,z)} \ldots \frac{I(G\sm \set{u, u_1 \dd u_k},z)}{I(G\sm \set{u, u_1 \dd u_{k-1}} ,z)}.\]
By $k$ applications of \refeq{recur}, one can recursively construct
functions $g_j(z)$, $j=1 \dd \ell$ such that
     \beql{fuzrec}
     f_u(z) = \frac{z}{(1-z)^{d(v)-\ell} \prod_{j=1}^\ell (1- g_j(z))},
     \eeql
where $g_j(z)$ is not identity for all $j$. Again consider $S_n$ with $u$ as one of the leaves then
the function given in \refeq{fuex} can be re-written as
     \[f_u(z) = \frac{z}{1- \frac{z}{(1-z)^{n-1}}}.\]
The \dt{depth of $f_u(z)$} is one more than the maximum depth of $g_j(z)$'s, where the base case
     \beql{base}
     f_u(z) = \frac{z}{(1-z)^\ell},
     \eeql
for $\ell \ge 1$, has depth one. The example function for $S_n$ with $u$ as a leaf
has depth two. The reason why we treat powers of $(1-z)$ in the denominator
separately will shortly become clearer.

Now if $\rho$ is a root of $I(G,z)$, then from \refPro{prop} we know that $f_u(\rho) =1$.
As mentioned in \refSec{intro}, to quantify the gap between $\beta(G)$ and the second
smallest absolute value over the remaining roots of $I(G,z)$, involves two steps.

Our first result is to show using \refPro{bloch} that $I(G,z)$ is injective in a
neighborhood of $\beta(G)$. For this purpose, we define 
        \beql{rbeta}
        r_G  \as \frac{\beta(G)^{\dia(G)}}{2n},
        \eeql
where $\dia(G)$ is the diameter of $G$ (see \refSec{prelim}). We show the following:

\bthml{mainone}
The polynomial $I(G,z)$ is injective on $D(\beta(G), r_G/2)$, that is, $\beta(G)$ 
is the unique root of the polynomial in this disc.
\ethml

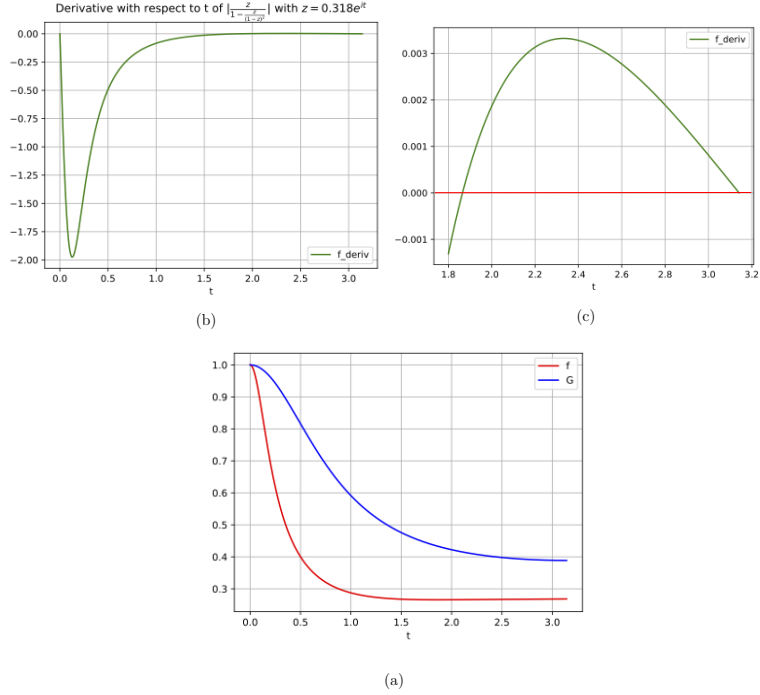
\vfigpdf{The absolute value function is not always monotone. In
(a) we have the plot of $\abs{z/(1-z/(1-z)^{2})}$ in red color, corresponding to $S_3$ as given in \refeq{fuex}, 
where $z\as \beta e^{i t}$, and $\beta\sim 0.318$ is the point where the function takes the value one. 
The plot in blue color shows the corresponding majorizing function.
(b) Shows the derivative with respect to $t$
of the absolute value in $[0, \pi]$ and (c) shows the same graph zoomed in to highlight 
an additional root of the derivative besides $0$ and $\pi$.}{example}{0.5}

We next need to show that for the points on the circle $\beta(G) \eitheta$ that are outside
the disc $D(\beta(G), r_G/4)$, there is a disc centered around each of the points such that
the value $\abs{f_u(z)}$ is smaller than one on these discs. The following result 
makes this precise:

\bleml{thetaa}
\label{lem:thetaa}
If $|w-\beta(G) \eitheta| \le r_G(1-\abs{f_{u}(\beta(G)\eitheta)})$ then $|f_u(w)|< 1$. 
In other words, there is no root of $I(G,z)$ in the disc 
$D(\beta(G) \eitheta, r_G (1-\abs{f_{u}(\beta(G)\eitheta)}))$, for $\theta > 0$.
\eleml

Note that the radius of the disc goes to zero as $\theta$ approaches zero, but this case is already
handled in \refThm{mainone}. In order to make the bound explicit in terms of graph parameters,
we need to upper bound $\abs{f_u(\beta\eitheta)}$ for $\theta$ sufficiently far away from the origin.
Ideally one would expect $\abs{f_u(\beta\eitheta)}$ to monotonically decrease in $\theta$ as it
varies from $0$ to $\pi$, but this is not the case, as shown in \refFig{example} for
the function $z/(1-z/(1-z)^{2})$.

Nevertheless, this is not far away from the truth, since as the depth of $f_u(z)$ increases
it concentrates around the origin and drops sharply as $\theta$ increases; we are not able
to prove this, but our observation is that it has the properties of  a ``good kernel'' \cite{stein-shakarachi}.
Instead, we show that there is a natural function that majorizes $\abs{f_u(z)}$ on the boundary
of the disc $D(0, r \eitheta)$, within its domain of holomorphy, and that is also
monotonically decreasing with $\theta$.

A \dt{majorant function} $G_r(\theta)$ for a complex valued function $g(r\eitheta)$ 
satisfies the following two properties:
\begin{enumerate}
\item $G_r(0)= g(r)$.
\item For all $\theta$, $|g(r\eitheta)| \le G_r(\theta)$, that is, the function majorizes
$g$ on the circle $r\eitheta$.
\item It is a monotone decreasing function, i.e., $G_r'(\theta) \le 0$; attains its maximum at the origin, i.e., $G'_r(0)=0$;
  and it is symmetric about the $y$-axis, i.e., $G_r(\theta) = G_r(-\theta)$.
\end{enumerate}
For example, consider the function in the base case $z/(1-z)^\ell$.
Its absolute value on $r\eitheta$ is 
$ r/|1-r\eitheta|^\ell$, which by a simple calculation turns out to
be 
     \[g_r(\theta) = \frac{r}{(1-2r\cos \theta + r^2)^{\ell/2}}.\]
In this case it is 
easier to argue monotonicity because the derivative with respect to $\theta$ is 
     \[\frac{-\ell r^2\sin\theta}{(1-2r\cos\theta+r^2)^{(\ell+2)/2}},\]
which is negative for $\theta \in (0, \pi)$. 
% Therefore, if we derive an upper bound on the function
% in a neighborhood of the origin, say $\theta \in [0, \theta_G)$, 
% that upper bound will hold for all $\theta > \theta_G$.
However, we use a simpler majorant
function that upper bounds the absolute value in the base case and behaves similarly. Based on the
observation that $\abs{1-r\eitheta} \ge 1 - r\cos \theta$, one such function is
      \beql{model}
      g_r(\theta) \as \frac{r}{(1-r\cos \theta)^\ell}.
      \eeql
For $r \le \beta(G-u)$, consider the function
     \[f_u(r\eitheta) = \frac{r\eitheta}{(1-r\eitheta)^\ell\prod_j (1-g_j(r\eitheta))},\] 
where $j$ varies over some fixed index, and $\ell \ge 0$. 
A {majorant function} $F_{u,r}(\theta)$ for $|f_u(r\eitheta)|$ is obtained recursively
from the majorant functions 
$G_{j,r}(\theta)$ for $\abs{g_j(r\eitheta)}$, respectively, as follows:
       \beql{model2}
       F_{u,r}(\theta) \as \frac{r}{(1-r\cos \theta)^\ell\prod_j (1-G_{j,r}(\theta))}.
       \eeql
The reason we treat powers of $(1-r\eitheta)$ separately is because if we take the absolute
value inside, as we will immanently do for the $g_j$'s, we will get a constant function $r$.
So, in principle, we assume that the $g_j$'s have {\em depth} more than one.

Let us verify that $F_{u,r}(\theta)$ satisfies all the properties of a majorant function.
Firstly, 
       \[F_{u,r}(0) 
       = \frac{r}{(1-r)^k\prod_j (1-G_{j,r}(0))}.\]
But as $G_{j,r}(0)=g_j(r)$ and the latter is positive it follows that
       \[F_{u,r}(0) = \frac{r}{(1-r)^k\prod_j |1-g_j(r)|}= f_u(r).\]
Secondly
     \[|f_u(r\eitheta)| 
       \le \frac{r}{|1-r\eitheta|^k\prod_j (1-\abs{g_j(r\eitheta)})}
       \le \frac{r}{(1-r\cos \theta)^k\prod_j (1-G_{j,r}(\theta))}
       = F_{u,r}(\theta)\]
and thirdly, taking the logarithmic derivative with respect to $\theta$ of $F_{u,r}(\theta)$ we obtain
      \beql{logFderiv}
       \frac{F_{u,r}'(\theta)}{F_{u,r}(\theta)} 
       = -\frac{kr \sin \theta}{(1-r\cos \theta)} + \sum_j \frac{G'_{j,r}(\theta)}{(1-G_{j,r}(\theta))}.
      \eeql
Therefore, 
the derivative is non-positive since by induction $G'_{j,r}(\theta) \le 0$ and
$F_{u,r}(\theta) > 0$; moreover, it also vanishes at $\theta = 0$, hence, the maximum
value is attained at the origin which is equal to $f_u(r) \le 1$, for $r \le \beta(G)$; the symmetric
nature also follows by applying induction to \refeq{model2}.

The majorant function will be used in \refLem{thetaa} instead of $\abs{f_u(\beta \eitheta)}$.
However, we still need a more explicit upper bound on $F_{u,r}(\theta)$. The monotone nature of
the function comes to rescue, since locally around the origin we will show that
the function is upper bounded by an inverted parabola, that is, 
$F_{u,r}(\theta) \le F_{u,r}(0) - c \theta^2$ for $\theta \le \theta_G$ 
and  some constant $c$ dependent on graph parameters.
Therefore, the upper bound at $\theta_G$ holds for all $F_{u,r}(\theta)$, for $\theta > \theta_G$.
Substituting this upper bound at $\theta_G$ in \refLem{thetaa} then gives us the desired explicit
disc around every point $\beta \eitheta$ that is devoid of roots.
In order to derive this local upper bound, we need to derive an upper bound on a variant
of the gamma-function for $F_{u,r}(\theta)$. This is done inductively. 
Define $F_u(\theta) \as F_{u, \beta(G)}(\theta)$, that is, $\fur(\theta)$ with $r=\beta(G)$, and
$G_j(\theta) \as G_{j,\beta(G)}(\theta)$. Then from \refeq{model2} it follows that
      \beql{futheta}
      F_u(\theta) = \frac{\beta(G)}{(1-\beta(G) \cos\theta)^\ell \prod_{j=1}^{d(u)-\ell} (1-G_j(\theta))}.
      \eeql

In particular, we
show the following result:
\bleml{inductivebound}
For $F_u(\theta)$ defined as in \refeq{futheta}, and 
   \[\Gamma \as \max_{j=1}^d \sup_{m \ge 0} \abs{\frac{G_j^{(m+1)}(0)}{(m+1)!}}^{1/(m+1)}\]
we have
     \[\sup_{k \ge 0} \abs{\frac{F_u^{(k)}(0)}{k!}}^{1/k} \le \frac{2d\Gamma}{\beta(G)}.\]
\eleml

As a consequence, we have 
\bleml{generalub}
If $\theta \le (\beta(G)/2d)^{2\Delta}$, where $d$ is the maximum degree of $G$,
and $\Delta$ is the depth of $F_u(\theta)$, then
   \[F_u(\theta) \le 1 - \frac{(\beta(G)\theta)^2}{4}.\]
To express the bound only in terms of $d$, we can use Shearer's bound $\lambda_S(d)$ 
instead of $\beta(G)$.
\eleml

Substituting  this in \refLem{thetaa}, we obtain 
\bcorl{gub}
Define
   \beql{thetaG}
   \theta_G \as \paren{\frac{\beta(G)}{4n}}^{\dia(G)}.
   \eeql
Then for all $\theta \ge \theta_G$
there is no root of $I(G,z)$ in the disc $D(\beta(G) \eitheta, r_G (\beta(G)\theta)^2/4)$,
\ecorl
Finally, combining this result with \refThm{mainone} gives us \refThm{main} as desired.

In the next sections, we develop the proofs and details of the results above.

\section{Univalence of $I(G,z)$ around $\beta(G)$}
\label{sec:univ-beta}
Throughout this section, we use $\beta \as \beta(G)$, and $v \in V$ as a representative
vertex. The function $f_v(z)$ given in \refeq{fuz} can also be expressed as
     \beql{fvz}
     f_v(z) = \frac{z I(G\sm N[v],z)}{I(G\sm v ,z)} = 1 - \frac{I(G,z)}{I(G\sm v, z)}.
     \eeql
Then $f_v(\beta)=1$, $f_v(0)=0$, and $f'_v(0)=1$. The next result gives a lower
bound on the growth of $f_v$ in the vicinity of $\beta_G$.

\blem 
At the smallest root $\beta(G)$ of $I(G,z)$, we have $f'_v(\beta_G) > 1/\beta_G$.
\elem
\bpf 
Taking the derivative on both sides of \refeq{fvz}, considering the second formulation, we get
that
        \[f_v'(z) = - \paren{\frac{I(G,z)}{I(G\sm v, z)}}' = -\frac{I'(G,z)}{I(G\sm v, z)} + \frac{I(G,z) I'(G\sm v,z)}{I(G\sm v, z)^2}.\]
Since $\beta$ is a root of $I(G)$, we get
        \[f_v'(\beta) = -\frac{I'(G,\beta)}{I(G\sm v,\beta)}.\]
Substituting \refeq{deriv} for the derivative, we further obtain that
\beql{fvbeta}
\begin{aligned}
          f_v'(\beta) &=  \sum_{u\in V}\frac{ I(G-N[u], \beta)}{I(G-v,\beta)}\\
                           &=  \frac{I(G-N[v], \beta)}{I(G-v, \beta)} + \sum_{u\in V, u \neq v}\frac{I(G-N[u], \beta)}{I(G-v,\beta)}\\
                           &=  \frac{1}{\beta} + \sum_{u\in V, u \neq v}\frac{I(G-N[u], \beta)}{I(G-v,\beta)},
\end{aligned}
\eeql
where the last step follows from the definition of $\beta$. Now observe that the graphs $G\sm N[u]$, for $u \in V\sm \set{v}$ 
and $G\sm v$ are all subgraphs of $G$. Therefore, from \refPro{prop} we know that 
their smallest root is larger than $\beta$ and hence 
the corresponding independence polynomials evaluated at $\beta$ are all positive. This means that 
the terms in the summation above are all positive, which gives us the desired lower bound on $f'_v(\beta)$.
\epf

We also have a corresponding upper bound on $f'_v(\beta)$.
\blem 
For all $v\in V$,  $f'_v(\beta) \leq \frac{n}{\beta^{\dia(G)}}$, where $\dia(G)$ is the diameter of $G$.
\elem
\bpf
Let the vertex set of the graph $G$ be $V={v,v_1, \ldots, v_{n-1}}$. 
% It is proven by [Csikvari], that the smallest root of the independence polynomial of a subgraph is 
% greater than the smallest root of the parent graph.
Again consider \refeq{fvbeta}. We start with a lower bound on the denominator $I(G-v,\beta)$ as follows.
We know that $I(G-v,\beta)=\beta I(G-N[v],\beta)$.
Let $G_v=G-N[v]+v_1$ where $v_1\in N[v]$, and $\beta_v$ be the smallest root of $I(G_v, z)$. 
Now $I(G,z)$ is continuously decreasing on the real line starting
from the origin down to its smallest root. At the origin we have 
$I'(G,0)=\abs{V(G)}$. Since $G_v$ is a subgraph of $G$, we have $\beta_v>\beta$. Therefore,
        \[I(G-v,\beta)=\beta I(G-N[v],\beta)=\beta I(G_v-v_1,\beta)\geq \beta I(G_v-v_1,\beta_v).\]
Repeating the above argument for vertex $v_1$ and so on we obtain $I(G-v,\beta)\geq \beta^k$,
where $k$ is the maximum distance of any vertex from the vertex $v$. Since the diameter $dia(G)$ 
of the graph is the longest shortest path in the graph and $\beta<1$ we have
        \beql{igvbeta}
        I(G-v,\beta)\geq \beta^{\dia(G)}.
        \eeql
Note that $I(G-N[u], \beta) < 1$ for all $u \in V$.
Substituting these two bounds  in \refeq{fvbeta} we get
        \[f'_v(\beta)\leq \frac{1}{\beta}+ \frac{n-1}{\beta^{\dia(G)}} 
                            = \frac{\beta^{\dia(G)-1}+ n-1}{\beta^{\dia(G)}}
                            \leq \frac{n}{\beta^{\dia(G)}},\]
since $\beta \leq 1$, $\beta^{\dia(G)} \leq 1$.
\epf

We next derive an upper bound on higher-order derivatives of $I(G,z)$, which will
be useful later.

\bleml{derivub}
Let $\beta$ be the smallest root of $I(G,z)$, then for all $0 \le k \le n$ 
we have $\abs{I^{(k)}(G, \beta)} \le {n \choose k} \le n^k$. 
More generally, if $H$ is a subgraph of $G$, then
$I^{(k)}(H, \beta) \le {|H| \choose k} \le |H|^k$.
\eleml
\bpf
% For $k = 0$ the claim is trivially true since left side is zero and right side is one. Let's understand
% the claim for $k=1$. In that case, from \refeq{deriv} we know that
%         \[I'(G, \beta) = - \sum_{u \in V} I(G-N[u], \beta).\]
% Now $G-N[u]$ is a subgraph of $G$ and hence by [Csivkari] its smallest root is larger than $\beta$.
% Since the independence polynomial goes down from one at the origin to zero at its smallest positive root
% it must be the case that $I(H, \beta) \le 1$, for all subgraphs $H$ of $G$. We will repeatedly use this 
% observation. In particular, each of the terms on the right-side above are smaller than one. Therefore, 
% we have the desired inequality for $k=1$.

% For $k=2$, we have 
%         \[I''(G,z) = \sum_{u \in V} I'(G-N[u]) = \sum_{u \in V} \sum_{v \in G-N[u]}I(G-N[u]-N[v], z).\]
% Again evaluating it at $\beta$, and using the basic observation above, we see that the sum on the right-side
% is at most $n^2$ (to be precise it is $\sum_u (n-d_u) = n^2 - 2|E|$).

% A similar argument shows that 
For a general $k$, we have
        \[I^{(k)}(G,z) = (-1)^k \sum_{\substack{u_1 \in V, u_2 \in G-N[u_1] \dd\\ u_k \in G-{\su_{j=1}^{k-1}}N[u_j]}} 
                                        I(G-N[u_1]-N[u_2]- \cdots N[u_k], z).\]
Since the graph $G-N[u_1]-N[u_2]- \cdots N[u_k]$ is a subgraph of $G$,  its evaluation at $\beta$ is smaller than one.
The number of distinct choices  of $u_1 \dd u_k$ are at most ${n \choose k}$, which completes
the proof.
\epf

Using the bounds above, we derive an upper bound on $\gamma_{I'(G,z)}$ (see \refeq{gamma}), the
standard gamma-function for the derivative $I'(G,z)$ at $\beta$.
We start with deriving a lower bound on $I'(G,\beta)$: Since
        \[I'(G, \beta) = - \sum_{u \in V} I(G-N[u], \beta)\]
and each $I(G-N[u], \beta)$ has  the same sign and by \refeq{igvbeta} is at least $\beta^{\dia(G)}$, we obtain
        \beql{iderivbeta}
        \abs{I'(G, \beta)} \ge n \beta^{\dia(G)}.      
        \eeql
It is not hard to see from the bound in \refLem{derivub} above and \refeq{iderivbeta} that
        \beql{gammab}
          \gamma =\gamma_{I'(G,z)}(\beta) \as \max_{k= 1\dd n}\abs{\frac{I^{(k+1)}(G,\beta)}{k!I'(G,\beta)}}^{1/k}
                                                 \le \max_{k= 1\dd n}\abs{\frac{n^{k+1}}{k! n \beta^{\dia(G)}}}^{1/k}
                                                  \le \frac{n}{\beta^{\dia(G)}}.
          \eeql

In order to apply \refPro{bloch} to $I(G,z)$ centered at $\beta$, we first need
to derive an upper bound on $\abs{I'(G,z)}$ in a neighborhood of $\beta$.
Consider the Taylor series expansion of the derivative around $\beta$
      \[I'(G,z) = \sum_{k \ge 0}\frac{I^{(k+1)}(G,\beta)}{k!} (z-\beta)^k.\]
Taking absolute values and pulling out the constant term we get
        \[
          \abs{I'(G,z)} \le \abs{I'(G,\beta)} \sum_{k \ge 0} \abs{\frac{I^{(k+1)}(G,\beta)}{k!I'(G,\beta)}} |z-\beta|^k.\]
Substituting the upper bound from \refeq{gammab} in the right-side, we obtain that for $z \in D(\beta, r)$
        \[
          \abs{I'(G,z)} \le \abs{I'(G,\beta)} \sum_{k \ge 0}\gamma^kr^k
                               = \frac{\abs{I'(G,\beta)}}{(1-r\gamma)},
                               \]
as long as $r\gamma < 1$. By \refeq{rbeta}, we know that, $r(G) \gamma \le 1/2$.
Therefore, for all $z \in D(\beta, r_G)$, $\abs{I'(G,z)} \le 2 \abs{I'(G,\beta)}$.
Substituting this upper bound, along with  the definition of $r(G)$,
in \refPro{bloch} applied to $I(G,z)$ at $\beta$, we get \refThm{mainone}.

\section{Gap around every point on the circle with radius $\beta(G)$}
\label{sec:gap-all-theta}
In this section, we prove \refLem{thetaa}, that is, 
we show that the gap to unity for every $\theta > \theta_G$ is governed
by the gap of $f_{u}(r\eitheta)$ to $f_{u}(r)$ and a constant that depends on $r$. We will
do the argument only for $r=\beta(G)$. For this purpose, we first need an upper bound
on  $\abs{f^{(k)}_u(\beta(G))}$. Again, for convenience, let $\beta \as \beta(G)$.

The $k$th derivative, up to sign, will have the form
        \beql{igzk}
        \paren{\frac{I(G,z)}{I(G-u,z)}}^{(k)} = \sum_{j=\max\set{0, k-n}}^k {k \choose j} I(G,z)^{(k-j)} \paren{\frac{1}{I(G-u,z)}}^{(j)}.
        \eeql
Applying Arbogast's formula \refeq{arbogast} we obtain 
        \beql{bderiv}
        \paren{\frac{1}{I(G-u,z)}}^{(j)}
        =\sum_{i_1 \dd i_j} \frac{(-1)^{i_1+ \cdots + i_j}}{I(G-u,z)^{1+i_1 + \cdots + i_j}}
        \frac{j!(i_1 + \cdots + i_j)!}{i_1! \ldots i_j!} 
        \prod_{m=1}^{j}\paren{\frac{I^{(m)}(G-u,z)}{m!}}^{i_m}
        % = \sum_{\ell=0}^j \frac{(-1)^\ell \ell!}{B(z)^{\ell+1}} \cdot
        %                 \calB_{j, \ell}(B'(z),B''(z)\dd B^{(j-\ell+1)}(z)),
        \eeql
% where $\calB_{j,\ell}(x_1 \dd x_{j-\ell+1})$ is the Bell polynomial
%         \beql{bell}
%         \calB_{j,\ell}(x_1 \dd x_{j-\ell+1}) \as \sum_{i_1 \dd i_{j-\ell+1}} \frac{\ell!}{i_1! \ldots i_{j-\ell+1}!} \prod_{m=1}^\ell \paren{\frac{x_m}{m!}}^{i_m}
%         \eeql
where the sum is over all  indices $i_1 \dd i_{j}$ such that 
        \beql{indicesj}
          i_1 + 2i_2 + 3 i_3 + \cdots + j i_{j} = j
        \eeql
Substituting $z=\beta$ in \refeq{igzk}, the term corresponding to $j=k$ disappears in the sum. 
Plugging the upper bound on the derivatives from \refLem{derivub} 
in \refeq{bderiv} above we get that
        \small\[
            \frac{\abs{f_u(\beta)^{(k)} }}{k!} \le \sum_{j=0}^{k-1} \frac{1}{j!(k- j)!} n^{(k-j)} 
            \cdot  \sum_{i_1 \dd i_j} \frac{1}{I(G-u, \beta)^{1+i_1 + \cdots + i_j}}\frac{j!(i_1 + \cdots + i_j)!}{i_1! \ldots i_j!}
                                                \prod_{m=1}^j \paren{\frac{n^m}{m!}}^{i_m}.
           \]\normalsize
From \refeq{indices}, we obtain that $n^{\sum_m mi_m} = n^j$. Moreover, as $I(G-u,\beta)<1$, 
we can upper bound its exponent by $k$ as well to get
        \[
            \frac{\abs{f_u(\beta)^{(k)} }}{k!}
                 \le \paren{\frac{n}{I(G-u, \beta)}}^k \;\sum_{j=0}^{k-1} \frac{1}{j!(k- j)!} 
                \sum_{i_1 \dd i_j}
                                        \frac{j!(i_1 + \cdots + i_j)!}{i_1! \ldots i_j!} \prod_{m=1}^j\paren{\frac{1}{m!}}^{i_m}.          \]
Using \refeq{partialBell} the summation over the indices $i_1 \dd i_j$ can be expressed as 
        \[
            \frac{\abs{f_u(\beta)^{(k)} }}{k!} 
            \le \paren{\frac{n}{I(G-u, \beta)}}^k\;\sum_{j=0}^{k-1} \frac{1}{j!(k- j)!} \sum_{t=0}^j t!B_{j,t}(1 \dd 1),
        \]
The last summation over $t$ is the ordered Bell number, $\wt{B}_j$ (see \refSec{prelim}), which gives us
        \[
            \frac{\abs{f_u(\beta)^{(k)} }}{k!} 
            \le \paren{\frac{n}{I(G-u, \beta)}}^k\;\sum_{j=0}^{k-1} \frac{\wt{B}_j}{j!(k- j)!}.
        \]
Furthermore, applying the upper bound from  \refeq{bellbound} we derive that
        \[
            \frac{\abs{f_u(\beta)^{(k)} }}{k!} 
            \le \paren{\frac{n}{I(G-u, \beta)}}^k\;\sum_{j=0}^{k-1} \frac{2^j}{(k- j)!}.
        \]
% using its exponential generating function $\sum_j \wt{B}_j/j! x^j = 1/(2-e^x)$ one gets 
% the bound that $\wt{B}_j/j! \le c^{j+1}$, for some constant $c > 1$; in fact, the precise statement is 
% that $\wt{B}_j/j! \sim \frac{1}{2(\ln 2)^{j+1}}$. Using the observation that for $c > 1$,
Notice that the summation
        \[\sum_{j=0}^{k-1} \frac{2^{j}}{(k-j)!}  
          = 2^k \sum_{j=1}^{k} \frac{1}{j!2^{j}}
          < 2^k \sum_{j=1}^{\infty} \frac{1}{j!2^{j}}
          = 2^k (e^{1/2}-1) \le 2^k.
        \]
Therefore, we finally obtain that
        \[
            \frac{\abs{f_u(\beta)^{(k)} }}{k!}
                 \le \paren{\frac{2n}{I(G-u, \beta)}}^k. \]
As a consequence, we get that
        \beql{gammafu}
        \gamma_{f_u}(\beta) \le \frac{2n}{I(G-u,\beta)} \le \frac{2n}{\beta^{\diaG}} = \frac{1}{r_G}.
        \eeql
From this bound, we can derive the following estimate, an alternate proof of \refThm{mainone}:
\blem
For all $z \in D(0, \beta+ r_G/2)$, $\abs{f_u(z)} \le 2$.
\elem
\bpf
A straightforward application of the triangle inequality to the Taylor series of $f_u$ around $\beta$ yields
        \[\abs{f_u(z)} \le f_u(\beta) \sum_{k \ge 0} (r\gamma_{f_u}(\beta))^k = \frac{f_u(\beta)}{1-r\gamma_{f_u}(\beta)}.\]
Since $r \le r_G/2$ it follows that $r\gamma_{f_u}(\beta) \le 2$, whence the upper bound
on $|f_u(z)|$ in $D(\beta, r_G/2)$. 
Because of the maximum modulus principle, the upper bound holds on
the whole disc. %$D(0, \beta+r_G/2)$.

\epf

Let $w$ be a point in the vicinity of $\beta(G) \eitheta$. Then taking the 
Taylor expansion around $\beta(G) \eitheta$, we get the following upper bound: 
       \[|f_{u}(w)| \le \abs{f_{u}(\beta\eitheta)} + \sum_{k \ge 1} \frac{|f^{(k)}_{u}(\beta \eitheta)|}{k!} |w-\beta \eitheta|^k.\]
Since $f_u$ is holomorphic with positive coefficients around the origin, from the strong maximum
modulus principle, the maximum of $|f^{(k)}_{u}(\beta \eitheta)|$ for all $\theta$ is attained
at the origin. Now, using the upper bound for $r=\beta$ from \refeq{gammafu}, we obtain that
       \[\abs{f_{u}(w)} < \abs{f_{u}(\beta\eitheta)} + \sum_{k \ge 1} {\left(\frac{|w-\beta \eitheta|}{r_G}\right)}^k.\]
Define $r \as |w-\beta \eitheta|/r_G$. If $r < 1$ then using the formula for a geometric series
we have
       \[\abs{f_{u}(w)} < \abs{f_{u}(\beta\eitheta)} + \frac{r}{1-r}.\]
Therefore, as long as 
       \[\frac{r}{1-r} \le 1- \abs{f_{u}(\beta\eitheta)}\]
we have $|f_u(w)| < 1$, or equivalently, if
       \[|w-\beta\eitheta| \le r_G \paren{\frac{1-\abs{f_{u}(\beta\eitheta)}}{2-\abs{f_{u}(\beta\eitheta)}}}\]
the function cannot take the value one in the vicinity of $\beta \eitheta$. The denominator can 
be simplified to one since the maximum value of $\abs{f_{u}(\beta\eitheta)}$ is one at the origin.
This completes the proof of the following \refLem{thetaa}.

% If $|w-\beta \eitheta| \le r_G(1-\abs{f_{u}(\beta\eitheta)})$ then $|f_u(w)|< 1$. In other words,
% there is no root of $I(G,z)$ in the disc $D(\beta \eitheta, r_G (1-\abs{f_{u}(\beta\eitheta)}))$,
% for $\theta > 0$.

\section{Upper bound on the Majorant Function near the origin}
\label{sec:escape}
In this section, we give the proofs of \refLem{inductivebound} and \refLem{generalub}.
We again use the shorthand $\beta \as \beta(G)$.
The idea is to consider the Taylor series expansion of $\fur(\theta)$ around the origin. 
More precisely,
we have
     \[\fur(\theta) = \fur(0) + \fur'(0) \theta + \sum_{k \ge 2} \frac{\fur^{(k)}(0)}{k!} \theta^k.\]
Since the first derivative vanishes, we have
     \[\fur(\theta) = \fur(0) + \sum_{k \ge 2} \frac{\fur^{(k)}(0)}{k!} \theta^k.\]
In fact, all the odd derivatives vanish and, the second derivative
is negative. We first verify the latter condition. From \refeq{logFderiv} it follows that
    \beql{secondderiv}
    \fur^{(2)}(\theta) = \fur'(\theta) \sum_j \frac{G'_{j,r}(\theta)}{(1-G_{j,r}(\theta))}
                           + \fur(\theta) \sum_j 
                                   \paren{\frac{G^{(2)}_{j,r}(\theta)}{(1-G_{j,r}(\theta))}
                                             + \frac{(G'_{j,r}(\theta))^2}{(1-G_{j,r}(\theta))^2}}.
    \eeql
Now inductively, the second derivatives $G_{j,r}^{(2)}(0)$ are negative (the base
case from \refeq{model} is $-r^2/(1-r)^2$), all the other terms vanish, 
which yield us that $\fur^{(2)}(\theta)(0)<0$ as desired.
This means that locally near the origin $\fur(\theta) \sim \fur(0) - c \theta^2$. We will next
show that for $\theta$ sufficiently small, half of the second term will dominate the
remaining summation in absolute value, and so $\fur(\theta) \le \fur(0) - \fur^{(2)}(0) \theta^2/4$
for $\theta \le \theta_G$.
For this purpose we need an upper bound on the absolute values of the $k$th derivatives 
of $\fur$ at the origin with $r=\beta$, which will be used to derive an upper bound on 
the $\gamma$-function for $\fur$. The upper bound will be derived inductively.

Let us begin with recalling some definitions: From \refeq{futheta} we have
     \[F_u(\theta) = \frac{\beta}{\prod_{j=1}^{d(u)} (1-G_j(\theta))},\]
where we have simplified the denominator to subsume the functions $\beta \cos\theta$
in the product by appropriate indexing, and $G_j(\theta) \as G_{j,\beta}(\theta)$.
We next derive a formula for the $k$th derivative of $F_u(\theta)$.

It can be verified that for $k \ge 1$, 
    \beql{fuderiv}
    F_u^{(k)}(\theta) 
    = \sum_{\ell=0}^{k-1} {k-1 \choose \ell} F_u^{(k-1-\ell)}(\theta) 
                 \sum_j - (\ln (1-G_j(\theta)))^{(\ell+1)}.
    \eeql
Using Arbogast's formula, \refeq{arbogast}, for the functions
$\ln(1-G_j(\theta))$, along with \refeq{logderiv}, we further get that 
    \small{\[
    F_u^{(k)}(\theta) 
    = \sum_{\ell=0}^{k-1} {k-1 \choose \ell} F_u^{(k-1-\ell)}(\theta) 
                 \sum_j \sum_{i_0 \dd i_\ell} 
                       \frac{(\ell+1)!}{i_0!\ldots i_\ell!}
                      \frac{(i_0 + \cdots + i_\ell)!}{(1-G_j(\theta))^{\sum_{m=0}^\ell i_m}} 
                      \prod_{m=0}^\ell \paren{\frac{G_j^{(m+1)}(\theta)}{(m+1)!}}^{i_m},
      \]}\normalsize
where $i_0 \dd i_\ell$ is an $(\ell+1)$-tuple of non-negative integers satisfying the equation
       \beql{indices2}
       i_0+ 2i_1 + \cdots + (\ell+1) i_\ell = \ell+1.
       \eeql
At this stage, we can inductively argue that if $k$ is odd then the derivative at the origin vanishes;
this is because one of the indices $m$ will be odd and by induction $G_j^{(m+1)}(0)$ vanishes.
Dividing both sides by $k!$, simplifying the binomial ${k-1 \choose \ell}$ term, 
and substituting $\theta=0$ we obtain that
    \beql{fuok}
    \frac{F_u^{(k)}(0)}{k!}
    = \sum_{\ell=0}^{k-1} \frac{F_u^{(k-1-\ell)}(0)}{k(k-1-\ell)! \ell!}
                 \sum_j \sum_{i_0 \dd i_\ell} 
                       \frac{(\ell+1)!}{i_0!\ldots i_\ell!}
                      \frac{(i_0 + \cdots + i_\ell)!}{(1-G_j(0))^{\sum_{m=0}^\ell i_m}} 
                      \prod_{m=0}^\ell \paren{\frac{G_j^{(m+1)}(0)}{(m+1)!}}^{i_m},
      \eeql
Define 
     \beql{gaamma}
     \Gamma \as \max_{j=1}^d \sup_{m \ge 0} \paren{\frac{\abs{G_j^{(m+1)}(0)}}{(m+1)!}}^{1/(m+1)},
     \eeql
which implies that 
     \[\frac{\abs{G_j^{(m+1)}(0)}}{(m+1)!} \le \Gamma^{m+1}.\]
Furthermore, we inductively assume that 
     \[\max_{j=0}^{k-1} \paren{\frac{\abs{F_u^{(j)}(0)}}{j!}}^{1/j} \le \frac{2d\Gamma}{\beta}.\]
% We will inductively assume that $F_u^{(j)}(0)/j! \le  (cdG/\beta)^j$, for some constant $c> 1$.
Since $F_u(0)=1$, we also know that 
$\prod_j (1-G_j(0)) = \beta$, which implies that for all $j$, $(1-G_j(0)) \ge \beta$.
Taking the absolute value, applying the triangle inequality, and
substituting these upper and lower bounds in the right-hand side, we get the following
\begin{align*}
    \frac{\abs{F_u^{(k)}(0)}}{k!}
    &\le \sum_{\ell=0}^{k-1} \paren{\frac{2d\Gamma}{\beta}}^{(k-1-\ell)}  \frac{1}{k\ell!}
                 \sum_j \sum_{i_0 \dd i_\ell} 
                       \frac{(i_0 + \cdots + i_\ell)!(\ell+1)!}{i_0!\ldots i_\ell!}
                      \prod_{m=0}^\ell 
                      \paren{\frac{\Gamma^{(m+1)}}{\beta}}^{i_m}
\end{align*}
From \refeq{indices2}, the term $\prod_m\Gamma^{i_m(m+1)} = \Gamma^{\ell+1}$ and hence
\begin{align*}
    \frac{\abs{F_u^{(k)}(0)}}{k!}
    &\le  \Gamma^{k}\sum_{\ell=0}^{k-1} \paren{\frac{2d}{\beta}}^{(k-1-\ell)}\frac{(\ell+1)!}{k\ell!}
                 \sum_j \sum_{i_0 \dd i_\ell} 
                       \frac{(i_0 + \cdots + i_\ell)!}{i_0!\ldots i_\ell!}
                     \paren{\frac{1}{\beta}}^{\sum_m i_m}\\
    &=  \Gamma^{k}\sum_{\ell=0}^{k-1} \paren{\frac{2d}{\beta}}^{(k-1-\ell)}\frac{(\ell+1)}{k}
                 \sum_j \sum_{i_0 \dd i_\ell} 
                       \frac{(i_0 + \cdots + i_\ell)!}{i_0!\ldots i_\ell!}
                      \prod_{m=0}^\ell 
                      \paren{\frac{1}{\beta}}^{\sum_m i_m}.
\end{align*}
The summation term over $j$ is independent of it, so we can upper bound the summation by the degree $d$.
Furthermore, if we define $t=i_0+ \cdots + i_\ell$, for a fixed $t$,
then the right-hand side further simplifies to
    \begin{align*}
    \frac{\abs{F_u^{(k)}(0)}}{k!}
    &\le d \Gamma^{k}\sum_{\ell=0}^{k-1} \paren{\frac{2d}{\beta}}^{(k-1-\ell)}\frac{(\ell+1)}{k}
                      \sum_{t=1}^{\ell+1}\frac{1}{\beta^t}
      \sum_{\substack{i_0 \dd i_\ell\\ i_0+\cdots+i_\ell=t}}
                       \frac{t!}{i_0!\ldots i_\ell!}\\
    &\le d\Gamma^{k}\sum_{\ell=0}^{k-1} \paren{\frac{2d}{\beta}}^{(k-1-\ell)}\frac{(\ell+1)}{k}
                      \sum_{t=1}^{\ell+1}\frac{1}{\beta^t} {\ell+1 \choose t-1},
     \end{align*}
where the last step follows from \refeq{bellbound3}. Since $\beta \le 1$, we can use $\beta^{-(\ell+1)}$ instead of $\beta^{-t}$
to  get
    \begin{align*}
    \frac{\abs{F_u^{(k)}(0)}}{k!}
    &\le d\paren{\frac{\Gamma}{\beta}}^{k}\sum_{\ell=0}^{k-1} \paren{2d}^{(k-1-\ell)}\frac{(\ell+1)}{k}
                      \sum_{t=1}^{\ell+1}{\ell +1\choose t-1}\\
    &\le d\paren{\frac{2\Gamma}{\beta}}^{k}\sum_{\ell=0}^{k-1} \frac{(\ell+1)}{k}d^{(k-1-\ell)},
     \end{align*}
where in the last step we upper bound the summation of the binomials by $2^{\ell+1}$.
Since $d \ge 2$, it can be showed that the new summation term is at most $d^{k-1}$, which
finally yields us the desired claim in \refLem{inductivebound}.

%\inductivebound*
    % \begin{align*}
    % \frac{\abs{F_u^{(k)}(0)}}{k!}
    % &\le  \paren{\frac{2d\Gamma}{\beta}}^k.
    %  \end{align*}
In the base case $\Gamma$
is the standard gamma-function for the cosine function, which is smaller than $1/\sqrt{2} < 1$.
Therefore, we get
    \begin{align*}
    \frac{\abs{F_u^{(k)}(0)}}{k!}^{1/k}
    &\le  \paren{\frac{2d}{\beta}}^{\Delta},
     \end{align*}
where $\Delta$ is the depth of $f_u(z)$.

In order for half of the second term in the Taylor series expansion of $F_u(\theta)$ around the origin
to dominate the sum of the remaining terms,
we want
   \[\abs{\frac{F_u^{(2)}(0)}{4}}\theta^2 \ge \sum_{\ell \ge 2} \paren{\frac{2d}{\beta}}^{2\Delta\ell}\theta^{2\ell}.\]
Assuming 
     \beql{theta1}
     \theta^2 \le (\beta/2d)^\Delta/2,
     \eeql
the above inequality follows if
    \beql{theta2}
    \theta^2 \le \frac{\abs{F_u^{(2)}(0)}}{4} \paren{\frac{\beta}{2d}}^\Delta.
    \eeql
We next derive an explicit lower bound on the second-derivative.

Recall that the first derivatives $F_u'(0)$ and $G_j'(0)$ vanish, and that the second derivatives
$G_j^{(2)}(0)$ are all negative. Therefore, from \refeq{secondderiv} we obtain that
      \[
        F_u^{(2)}(0) = \sum_j \frac{G_j^{(2)}(0)}{(1-G_j(0))} 
                     \ge \sum_j G_j^{(2)}(0).\]
Therefore, $\abs{F_u^{(2)}(0)}\ge d \min \abs{G_j^{(2)}(0)}$. Inductively, we obtain that
      \[\abs{F_u^{(2)}(0)}\ge d^{\Delta-1}\beta^2/(1-\beta)^{d+1},\] 
since the absolute value of the second derivative of the 
majorant function in the base case is $d\beta^2/(1-\beta)^{d+1}$.
Substituting this in  \refeq{theta2}, we get a slightly weaker constraint than \refeq{theta1}, namely,
$\theta^2 \le {\beta^{\Delta+2}}/{4d}$.
So, in order to simplify, we combine the two constraints to obtain \refLem{generalub}.
%\generalub*

In order to combine this lemma with \refThm{mainone}, we notice that the disc $D(\beta, r_G/2)$
subtends the angle $\arcsin(r_G/2\beta)$, which is at least $\beta^{\dia(G)}/(4n)$, at the origin. 
Take
$\theta_G$ as a quantity smaller than this bound and the constraint in \refLem{generalub}, namely
as defined in \refeq{thetaG}.
Since the function $F_u(\theta)$ is monotonically decreasing, we know that
for all $\theta \ge \theta_G$, $F_u(\theta) \le F_u(\theta_G)$.
Substituting  this in \refLem{thetaa}, we obtain that for all $\theta \ge \theta_G$,
the disc 
     \[D(\beta \eitheta, r_G(\beta\theta_G)^2/4)\] 
is devoid of roots. Since the
radius here is smaller than $r_G/2$, we combine this with \refThm{mainone} to finally obtain that
the disc
   \[D\paren{0, \beta + \paren{ \frac{\beta}{n}}^{O(\dia(G))}}\]
contains exactly one root of $I(G,z)$ completing the proof of \refThm{main}. Note that
the depth $\Delta$ is smaller than the diameter $\dia(G)$, which instead is bounded by $n$.
%%
%% bibliography
%%

%% Please use bibtex, 
%\section{Examples of Root Gaps}
%In the case of cycle graph $C_n$, path graph $P_n$ and the complete bipartite graph $K_{n \times n}$ the gap between $\beta(G)$ and the absolute value of the next root is $1+ \Omega\left(\frac{1}{n}\right)$.
\section{Some Explicit Lower Bounds On The Gap}\label{sec:example}
In this section, we derive explicit lower bounds on the gap between
the smallest root $\beta(G)$ and the root with the second smallest
absolute value of the independence polynomial of some fundamental
graph classes. In particular, we give the explicit lower bounds for
the Path Graph $(P_n)$, the Cycle Graph $(C_n)$ on $n$ vertices and
the Complete Bipartite Graph $(K_{n \times n})$ on $2n$ vertices. For
convenience, let $\beta \as \beta(G)$ and $\alpha \as \alpha(G)$ be the root with the second
smallest absolute value of the independence polynomial of graph $G$ in each of the cases.

\paragraph*{Path Graph}
To describe the independence polynomial $P_n(z)$ of the Path Graph we need the
Fibonacci polynomials \cite{Hoggatt01101973}: Let $F_1(z) \as 1$, $F_2(z) \as z$
and recursively define $F_{n+1}(z) \as z F_n(z) + F_{n-1}(z)$.
From the relation \refeq{recur} for $P_n(z)$ we obtain
that $z^{n-1}P_n(-1/z^2)= F_{n+2}(z)$. The roots of $F_{n+2}(z)$ are $2i \cos(k\pi/(n+2))$,
$k=1 \dd n+1$. Therefore, the roots of $P_n(z)$ are $1/(4\cos^2(k\pi/(n+2)))$,
for $k=1 \dd \lfloor \frac{n+1}{2}\rfloor$. Therefore, $\beta = \frac{1}{4 \cos^2\left(
\frac{\pi}{m}\right)}$ and $\alpha = \frac{1}{4 \cos^2\left(
\frac{2\pi}{m}\right)}$, with $m\as n+2$.
% , and the ratio $\frac{\alpha}{\beta}=
% \frac{\cos^2\left( \frac{\pi}{m}\right)}{ \cos^2\left(
% \frac{2\pi}{m}\right)}$, with $m = n+2$. 
Using
the Taylor series for the cosine function, for large values of $m$, we get
% \begin{align*}
% \cos\left( \frac{\pi}{m}\right) 
% &= 1 - \frac{\pi^2}{2m^2} + O(m^{-4}) \\
% \mbox{and, }\; 
% \cos\left( \frac{2 \pi}{m}\right) &= 1 - \frac{2 \pi^2}{m^2} + O(m^{-4}) \\
% \end{align*}
   \[\frac{\cos\left( \frac{\pi}{m}\right)}{\cos\left( \frac{2 \pi}{m}\right)} = 1 + \frac{3\pi^2}{2m^2} + O(m^{-4}).\]
On squaring, we obtain that 
   \[\frac{\alpha}{\beta} 
     = \frac{\cos^2\left( \frac{\pi}{m}\right)}{\cos^2\left( \frac{2 \pi}{m}\right)}= 1 + \frac{3 \pi^2}{m^2} + O(m^{-4})
     = 1 + \Omega\left( \frac{1}{n^2}\right).\]

\paragraph*{Cycle Graph}
The independence polynomial $C_n(z)$ of the Cycle Graph can be expressed in terms of
the Chebyshev polynomial of the first kind $T_n(z)$ \cite{levit}:
    \[C_n(z) = 2 z^{n/2} T_n\paren{\frac{1}{2\sqrt{z}}}.\]
Therefore, its roots are $\frac{1}{4 \cos^2\left( \frac{(2k+1)\pi}{2n}\right)}$, 
$k= 0,1, \ldots, \lfloor \frac{n-1}{2}\rfloor$. 
This implies that $\beta =
\frac{1}{4\cos^2\left(\frac{\pi}{2n}\right)}$ and $\alpha =
\frac{1}{4\cos^2\left(\frac{3\pi}{2n}\right)}$. An argument similar to the one above
implies that asymptotically we have
% Let $a= \frac{\pi}{2n}
% $, again as before, asymptotic expansion using the Taylor series for
% $\cos$ function gives us
% \begin{align*}
% \cos(a) &= 1- \frac{a^2}{2} + \frac{a^4}{24} + O(a^6)\\
% \mbox{and, }\; \cos(3a)&= 1- \frac{9a^2}{2} + \frac{81a^4}{24} + O(a^6)\\
% \mbox{Thus, }\; \frac{\cos(a)}{\cos(3a)} &= 1+ 4a^2 + \frac{44}{3}a^4 + O(a^6)\\
% \end{align*}
   \[\frac{\alpha}{\beta} 
     = 1+ \frac{2\pi^2}{n^2} + \frac{17\pi^4}{6 n^4} + O(n^{-6}) 
     = 1 + \Omega\paren{\frac{1}{n^2}}.\]

\paragraph*{Complete Bipartite Graph}
The independence polynomial $K_{n\times n}(z)$ of the Complete Bipartite Graph on $2n$
vertices is $K_{n \times n}(z) = 2{(1-z)}^n-1$. 
Its roots are $z_k = 1- 2^{-\frac{1}{n}}e^{\frac{i2k\pi}{n}}$, $k=0,1,\ldots, n-1$. 
Therefore $\beta = 1- 2^{-\frac{1}{n}}$ and
$\alpha = 1-re^{i\frac{2\pi}{n}}$, where $r \as 2^{-\frac{1}{n}}$. 
Since, in this case $\alpha$ is a complex root, we compute ratio of the absolute values of
the roots:
\begin{align*} 
  \frac{\left|\alpha\right|}{\beta} = \frac{\sqrt{1+r^2- 2r\cos(\frac{2\pi}{n})}}{1-r}
                                   =\frac{\sqrt{{(1-r)}^2 + 4r\sin^2(\theta/2)}}{1-r} 
                                   =\sqrt{1+\frac{4r}{({1-r}^2)}\sin^2\left(\theta/2\right)}.
\end{align*} 
Take $b \as \ln(2)$ then, $r = e^{-b/n} = 1- \frac{b}{n} +
\frac{b^2}{2n^2} + O(n^{-3})$ and $\sin(\frac{\pi}{n})= \frac{\pi}{n}
- \frac{\pi^3}{6n^3} + O(n^{-5})$. After substituting these estimates and
simplifying we get that for large $n$
\begin{align*} 
\frac{\left|\alpha\right|}{\beta}
  & = \sqrt{1 +\frac{4\pi^2}{b^2} + \left(\frac{-\pi^2}{3} - \frac{4\pi^4}{3b^2}\right)\frac{1}{n^2} + O(n^{-3})}\\ 
  &= \sqrt{1 + \frac{4\pi^2}{b^2}} -\frac{\pi^2}{6}\frac{1}{n^2} + O(n^{-3}) \\
  &\approx 9.119 - O\left(\frac{1}{n^2}\right).
\end{align*}

\section{Concluding Remarks}
This paper provides the first quantitative lower bound on the gap between the smallest
root of the independence polynomial and the second smallest absolute value. A simple proof for the existence of the gap can be based on the maximum modulus principle.
Quantifying the principle for the special function at hand involves studying its local
behaviour. To the best of our knowledge, this is the first time a result is provided for separation
between roots of the independence polynomial; earlier results always focus on zero-free regions. Our 
larger hope is to study the algorithmic implications of the ratio of $\beta(G)$ to the
second smallest absolute value. Can it be used to design algorithms that are efficient for
those graphs where this ratio is large, for example, the graph classes mentioned in \refSec{example}?

\bibliographystyle{abbrv}
\bibliography{local}
% \appendix

% \section{}\label{sec:itemStyles}

\appendix
% \newpage
\section{Basic Results}
\label{sec:prelim}
Throughout this paper $G=(V,E)$ will be assumed to be a simple undirected graph.
Let $n$ be the number of vertices in $G$. For every vertex $u \in V$, 
let $d(u)$ denote its degree, and $d\as \max_u d(u)$, for all $u$; we will assume that $d \ge 2$. 
Let $\dia(G)$ denote the \dt{diameter of $G$}, that is the length of the longest path between any 
pair of vertices in any connected component of $G$. In the subsequent sections, we will need some basic
properties of the independence polynomial (see \cite{csikvari:note:12,bencs+3:23} for proofs):

\bpropl{prop}
Let $G$ be as above and $I(G,z)$ be its independence polynomial.
\begin{enumerate}
% \item If $G$ is disconnected with components $H_1 \dd H_k$ then
%        \[I(G,z) = \prod_{j=1}^k I(H_i,z).\]
\item The derivative $I'(G,z)$ satisfies
        \beql{deriv}
        I'(G,z) = - \sum_{u \in V}I(G-N[u], z).
        \eeql
\item If $H$ is a subgraph of $G$ then $\beta(G) \le \beta(H)$.
\item Shearer \cite{shearer:85} showed a lower bound on $\beta(G)$, namely,
    \beql{shearer}
    \beta(G) \ge \lambda_S(d) \as \frac{(d-1)^{d-1}}{d^d}.
    \eeql
     Applying the third property with $G$ as the complete graph and $H$ an arbitrary graph
     with at most $n$ vertices, we also have the following lower bound 
     \beql{betalb}
     \beta(H) \ge \frac{1}{n}.
     \eeql
\end{enumerate}

\epropl

We will also need variants of Smale's gamma-function \cite{bcss:bk}: 
Given a function $f: \CC \mt \CC$ holomorphic
at a point $z \in \CC$, such that $f^{(j)}(z) \neq 0$, define
     \beql{gamma}
     \gamma_{f,j}(z) \as \sup_{k > j} \abs{\frac{j!f^{(k)}(z)}{k!f^{(j)}(z)}}^{1/(k-j)}.     
     \eeql
The standard gamma-function corresponds to $\gamma_{f,0}(z)$, and we will simply use
$\gamma_f(z)$ to denote that. Intuitively, the function is related to the inverse of the
radius of convergence of $f^{(j)}(z)$ at $z$.

We recall Arbogast's formula (also called the formula of Fa\`a di Bruno) \cite{krantz-parks} 
for derivatives of composition of functions: For $N \in \ZZ_{\ge 0}$, the 
$N$th derivative of
        \beql{arbogast}
         (f \circ g)^{(N)}(z) 
        = \sum_{i_1 \dd i_N} \frac{N!}{i_1! \ldots i_N!} f^{(i_1 + \cdots + i_N)}(g(z)) \prod_{m=1}^N \paren{\frac{g^{(m)}(z)}{m!}}^{i_m},
        \eeql
where the sum is over all tuples of non-negative integers $i_1 \dd i_N$ such that
        \[
        i_1 + 2i_2 + 3 i_3 + \cdots + Ni_N = N.
        \]
Given a $K$, we can combine the terms corresponding to $i_1 + \cdots + i_N=K$  and simplify
the summation as follows. Since $\sum_j j i_j = N$ it follows with the additional 
constraint that $\sum_j(j-1)i_j = N-K$. This implies that 
$i_j = 0$, for $j > N-K+1$, and so we can express \refeq{arbogast} as
       \beql{arbo2}
        (f\circ g)^{(N)}(z) 
        = \sum_{K=0}^N \sum_{i_1 \dd i_N} f^{(K)}(g(z)) B_{N,K}(g'(z), g^{(2)}(z) \dd g^{(N-K+1)}(z)),
        \eeql
where $B_{N,K}(x_1 \dd x_{N-K+1})$ are \dt{the partial exponential Bell polynomials}:
        \beql{partialBell}
        B_{N,K}(x_1 \dd x_{N-K+1}) \as \sum_{i_1 \dd i_{N-K+1}} \frac{N!}{i_1! \ldots i_{N-K+1}!} \prod_{m=1}^{N-K+1} \paren{\frac{x_m}{m!}}^{i_m},
        \eeql
and $i_1 \dd i_{N-K+1}$ satisfy the following two constraints:
        \beql{indices}
        \begin{split}
          i_1 + 2i_2 + 3 i_3 + \cdots + (N-K+1) i_{N-K+1} = N \text{ and }\\
          i_1 + i_2 +  i_3 + \cdots + i_{N-K+1} = K.
        \end{split}
        \eeql
Notice that $B_{N,K}(1 \dd 1) = \stirtwo{N}{K}$, the Sterling number of second kind, 
that is number of ways to partition an $N$ element set into $K$ non-empty parts, and hence
        \[\sum_{K=0}^NK! B_{N,K}(1 \dd 1) = \wt{B}_N\] 
the \dt{ordered Bell number}, which satisfy the following recurrence:
        \beql{bellrecur}
        \wt{B}_N = \sum_{i=1}^{N-1} {N \choose i}\wt{B}_{i}.
        \eeql
From this, we can derive the following claim inductively: 
        \beql{bellbound}
        \frac{\wt{B}_N}{N!} \le \paren{\frac{1}{\ln 2}}^N.
        \eeql
% Let us suppose that for $i < N$ the claim is true with some $c >1$. 
% Then from the recurrence formula we obtain:
%         \[
%           \frac{\wt{B}_N}{N!} \le \sum_{i=1}^{N-1} \frac{c^i}{(N-i)!} 
%                               = c^N \sum_{i=1}^{N-1} \frac{1}{c^{N-i}(N-i)!}
%                               < c^N (e^{1/c}-1),
%         \]
% which is smaller than $c^N$ as long as $e^{1/c} \le 2$, that is, $c \ge \frac{1}{\ln 2} \sim 1.44$.
% We will also need a variant of \refeq{bellbound}: for $N \ge 1$,
%         \beql{bellbound2}
%         \frac{\wt{B}_N}{(N-1)!} \le 3^N.
%         \eeql
% This follows from the earlier bound and the fact that $3 \ln 2 > 2$, 
% since then $N \le (3 \ln 2)^N$, for $N \ge 1$.

We will also need the following observation:
       \beql{bellbound3}
       \sum_{i_1 \dd i_{N-K+1}} \frac{K!}{i_1! \ldots i_{N-K+1}!} = {N-1 \choose K-1}
       \eeql
where the sum is over all tuples $(i_1 \dd i_N)$ satisfying the conditions in \refeq{indices}.
The left-hand side counts all partition of $N$ into
$K$ blocks where the ordering of distinct blocks only matter (distinct blocks correspond 
to different choices of $j$). However, these are precisely the number of compositions of $N$
into $K$ parts, which is the term on the right-hand side.

Besides the above, we need the following observations for $j \ge 0$
%  on higher derivatives of basic functions, 
% which can be inductively verified:
        \beql{invderiv}
        \paren{\frac{1}{x}}^{(j)} = \frac{(-1)^j j!}{x^{j+1}}
        \eeql
and
        \beql{logderiv}
        \paren{\ln x}^{(j)} = \paren{\frac{1}{x}}^{(j-1)} =  \frac{(-1)^{j-1} (j-1)!}{x^{j}}.
        \eeql
Given $\alpha \in \CC$ and $r > 0$, we will denote by $D(\alpha, r) \ibp \CC$
the open disc centered at $\alpha$ with radius $r$.

In addition to the above, we need the following quantified result on 
the injectiveness of a function in a neighborhood of a point where the derivative does not vanish
\cite{Harris1977}:
\bpropl{bloch}
Let $h: D(\alpha, r) \mt \CC$ be a holomorphic function such that $h(\alpha)=0$ and $h'(\alpha) \neq 0$.
Then $h$ is injective on the disc 
\[D\paren{\alpha, \frac{r|h'(\alpha)|}{\sup_{z \in D(\alpha, r)}|h'(z)|}}.\]
\epropl

\end{document}